# Use of Honey-Bee Mating Optimization Algorithm to Design Water Distribution Network in Gurudeniya Service Zone, Gurudeniya, Sri Lanka

### K.H.M.R.N. Senavirathna, S. Thalagala and C.K. Walgampaya

**Abstract:** Water distribution network (WDN) is a highly complex urban infrastructure that is intended to supply water from source node to consumer nodes. As far as the components of an infrastructure are concerned, the interconnecting pipes that transport water from source node to the demand nodes account for the major fraction of the capital cost. Therefore, when designing WDN, the primary concern is on finding the optimal combination of pipe diameters that satisfies the hydraulic-head requirements at minimum cost. More recent methodologies employ direct heuristic and stochastic heuristic algorithms for the design of WDN. In this study, the use of Honey-Bee Mating Optimization (HBMO) stochastic algorithm for the design of WDN in the Gurudeniya Service Zone, Gurudeniya, Sri Lanka, has been examined. In-built hydraulic constraints have enabled the algorithm to deliver results without needing a separate hydraulic simulation software in this design. The solutions given by HBMO algorithm and the solution implemented by National Water Supply and Drainage Board, Sri Lanka have been compared. From the implementation of the HBMO algorithm, results show that the HBMO algorithm is successful in addressing the WDN design problem with given hydraulic constraints for Gurudeniya Service Zone.

**Keywords:** Honey-bee mating optimization, Water distribution network, Water resources management

## 1. Introduction

Water distribution network is an infrastructure that is intended to collect, transmit, treat, store, and distribute water for homes, commercial establishments, industry, and irrigation, as well as for such public needs as firefighting and street flushing. Of all municipal services, provision of potable water is considered to be the most vital department. People depend on a water supply system for drinking, cooking, washing, carrying away wastes, and other domestic needs. Water supply systems shall comply with the requirements for public, commercial, and industrial activities. In all cases, the water delivered through the system must fulfill both quality and quantity requirements [12].

Water utilities seek to provide customers with a safe, reliable, continuous supply of high quality water while minimizing costs. Water is often delivered through a complex distribution system involving miles of pipe and often incorporating numerous pumps, regulating valves, consumer nodes and storage reservoirs. The performance of these systems is often difficult to be understood due to its physical size and complexity and presence of large amount of system information and data when determining its function.

Therefore, analyzing and modeling of these water distribution systems are essential to meet the goals of delivering safe, reliable water supply at optimal cost. Network analysis is critical for proper operation and maintenance of a water supply system. [3].

In practice, the conventional way of designing a water distribution network is manual analysis by making an educated guess for the least cost combination of pipe diameters.

*Eng. K.H.M.R.N. Senavirathna, B. Sc. Eng (Hons) (Peradeniya), Temporary Instructor, Department of Engineering Mathematics, Faculty of Engineering, University of Peradeniya.*
*Email:rajithas@eng.pdn.ac.lk*
*ORCID ID: https://orcid.org/0000-0002-4081-3781*
*Eng. S. Thalagala, B. Sc. Eng (Hons) (Peradeniya), Temporary Instructor, Computing Center, Faculty of Engineering, University of Peradeniya.*
*Email:shiron.thalagala@eng.pdn.ac.lk*
*ORCID ID: https://orcid.org/0000-0002-0396-1849*
*Eng. (Dr.) C.K. Walgampaya, AMIE (Sri Lanka), B.Sc. Eng. (Hons) (Peradeniya), M.Sc.(UofL,KY, USA), Ph.D. (UofL, KY, USA), Senior Lecturer, Department of Engineering Mathematics, Faculty of Engineering, University of Peradeniya.*
*Email:ckw@pdn.ac.lk*
*ORCID ID: https://orcid.org/0000-0003-2192-474X*

In the recent past, the rise of water simulation software, such as WaterGEMS, WaterCAD, and EPANET, has enabled WDN design engineers to design water distribution networks with relatively less effort. However, a globally accepted right way of determining optimal combination of pipe diameters, especially for a large scale water distribution network, is still a subject of research due in part to the complexity of the problem. A small change in the combination of pipe diameters will result in a large change in the capital that has to be invested in the project implementation [11].

## 2. Literature Review

Due to these issues, instead of making educated guesses for optimal combination of pipe diameters, notable amount of research shows that WDN design engineers and analysts have started employing direct-heuristic and stochastic-heuristic algorithms for the design of water distribution networks.

Noticeable studies that have used direct methods include the works done by Alporevits and Shamir [1], Suribabu [10], Lin et al. [7], Hsu and Cheng [6], Mohan and Babu [8], and Awe et al. [2]. As far as the general direct method criterion to determine pipe diameters of WDN is concerned, initially, all the pipes in the network have been set to their minimum discrete commercially available size, and the network has been simulated to determine the pressures at all the nodes in the network. The hydraulic simulations, alongside with the proposed algorithm, have been carried out with a hydraulic simulation software. If the minimum pressure requirements have not been satisfied in a node, the pipe exhibiting the maximum velocity has been increased to next commercially available size. Then the network has again been simulated to determine the pressure and flow velocity. Identification of maximum flow velocity in a pipe and increasing its size has been continued until the pressure head at all the demand nodes has exceeded the corresponding minimum required. At the end of this process, a feasible solution to the network has been obtained. The pipes for which the diameter has been more than required should have been replaced with smaller diameters to achieve the cost effective design. Hence, at the second stage, each pipe diameter has then been reduced to the next available size based on the selection indices proposed according to the study, until the pressure head at any one node has approached the minimum required. It is noticeable that, these general direct methods have been applicable where the shape of the objective function has been known with respect to the set of decision variables (i.e. diameters). On the other hand, as a drawback, if the shape of the objective function is unknown, these direct methods become partly misleading when aiming for a global optimal solution.

On the contrary to the direct methods, the stochastic methods employ a search procedure over the entire solution space. When picking a solution whose objective function to be evaluated, the stochastic algorithms use probabilistic rules inspired mostly by the biological or any other natural phenomena. The applicability even with the complex non-linear problems, not-needing to identify the shape of the objective function (i.e. convex, concave etc.), are some of the benefits of these stochastic-heuristic algorithms. A few recent studies that use stochastic algorithms such as, Simulated Annealing Approach by Cunha and Sousa [4], Genetic Algorithm (GA) by Dijk et al. [5], and Honey-Bee Mating Optimization (HBMO) by Mohan and Babu [9], have been carried out to check their performance upon given WDNs. Further, in the study done by Mohan and Babu [9], it has been systematically proven the efficacy of the HBMO algorithm over the other well-established stochastic algorithms by making use of the number of iterations taken to reach the optimal solution as a basis of evaluation.

In this study, the use HBMO algorithm for the design of WDN of Gurudeniya Service Zone, Gurudeniya, Sri Lanka, has been examined. The basic structure of the HBMO algorithm was adopted from Mohan and Babu [9]. To preserve the simplicity of the algorithm, it was updated with the change of the omission of the Nurse bees' mutation process. Further, comparison of solutions given by HBMO algorithm to the solution implemented by National Water Supply and Drainage Board (NWSDB), Sri Lanka, has been examined.

In the following sections, methodology and analysis, followed by the results and discussion, are presented. The conclusions made in the study are given towards the end of the paper.

## 3. Methodology and Analysis

In this section, the terminology in the study, mathematical model formulation, and adaptation of HBMO algorithm for the WDN

design of Gurudeniya Service Zone are written explicitly. A description of the type of data used in the study is given at the latter part of this section.

## 3.1 Terminology

- Water pipe:

A water pipe is a component for conveying water from source node to consumer demand nodes. This can be in meter scale or in kilometer scale.

- Consumer node or demand node:

This is a point from which water is taken out of the water distribution network. These nodes could be housing units, factories, or any other establishment that draws water from the water distribution network.

- Serial/looped pipe configuration:

A serial pipe configuration is a WDN configuration where all the pipes in the network are connected in series, whereas a looped configuration has at least one loop of pipes connected.

- Water source or reservoir or tank:

This is a source point by which water supplier distributes water through the network to the consumer nodes. Usually, this source point is chosen to be at a higher elevated place. An advantage of the source point being at a high elevated place is that the water can then be sent through the distribution network by the use of earth's gravitational pull. This kind of flow is commonly known as gravity flow. Pumps may be used when the water pressure at consumer nodes is not significant with the only use of gravity flow.

## 3.2 Optimization Model Formulation

The objective function of the design of water distribution network is defined as,

$$Min\ Z = \sum_{i=1}^{N} C_i(D, L) \quad \ldots(1)$$

where,
$N$ = number of pipes in the water distribution network;
$C_i(D, L)$ = cost of the i$^{th}$ pipe having diameter $D$ and length $L$.

The hydraulic constraints to be satisfied while minimizing the objective function are,

$$Q_i = W_j + Q_{i-1} \quad i = 1,2,3,\ldots,np \quad \ldots(2)$$
$$j = 1,2,3,\ldots,nd$$

$$H_i = H_{i-1} - F_i \quad i = 1,2,3,\ldots,np \quad \ldots(3)$$

$$H_{R\ j} = H_j - E_j \quad j = 1,2,3,\ldots,nd \quad \ldots(4)$$

$$H_{R\ j} \geq H_{R\ j}^{min} \quad j = 1,2,3,\ldots,nd \quad \ldots(5)$$

$$h_{f\ i} \leq h_{f\ max} \quad i = 1,2,3,\ldots,np \quad \ldots(6)$$

$$D_i \geq 0 \quad i = 1,2,3,\ldots,np \quad \ldots(7)$$

where,
$Q_i$ = flow in the $i^{th}$ pipe;
$W_j$ = demand in the $j^{th}$ node;
$H_i$ = nodal water head at $i^{th}$ node;
$H_{R\ j}$ = residual water head available at the $j^{th}$ node;
$E_j$ = nodal elevation (from mean sea level MSL) at $j^{th}$ node;
$H_{R\ j}^{min}$ = minimum residual water head required at the $j^{th}$ node;
$nd$ = number of demand nodes; and
$np$ = number of pipes.

The friction of a pipe is calculated using the Hazzen's Williams Formula as follows.

$$h_{f\ i} = C_{ft}\ \frac{10.666 * Q_i^{1.85}}{C_{HW}^{1.85} * D_i^{4.87}} \quad i = 1,2,3,\ldots,np \quad \ldots(8)$$

$$F_i = L_i * hf_i \quad i = 1,2,3,\ldots,np \quad \ldots(9)$$

where,
$h_{f\ i}$ = friction loss in the $i^{th}$ pipe;
$C_{ft}$ = fitting loss coefficient;
$C_{HW}$ = Hazzen's Williams Coefficient
$F_i$ = total friction loss of the $i^{th}$ pipe; and
$L_i$ = length of the $i^{th}$ pipe.

## 3.3 Honey-Bee Mating Optimization (HBMO) Algorithm

### 3.3.1 Biological Basis of Honey-Bee Colony

A honey-bee colony is composed of queens, drones, and workers or nurses. Division of labor, individual and group level communications, and cooperative behavior are the strange features of the honey-bee colony. The number of queen bees in a colony may be either one (monogynous) or many (polygynous). Queens are the main reproductive sources of the honey-bee colony. Drones (the male bees) are fathers of the colony. The role of a drone is to inject its sperm into queen's spermatheca. The nurse bees take care of broods' growth and do not involve in the process of reproduction. As a honey-bee colony is concerned, the mating process begins with the queen bee's mating flight and its travel far away from the nest. At the beginning of the mating flight, the queen's speed is very fast and the speed gradually gets reduced as the mating process is progressing. During the mating flight of the queen, swarm of drones follow the queen bee and mate with the queen in the air. The queen bee can have several number of matings with different drones even in a single mating flight. But a drone's life ends up with a single mating. Eventual death of the drone indicates the successful insemination. The sperms accumulated in the queen's spermatheca create a genetic pool for the new colony to be born. The queen starts off a mating flight with an empty spermatheca, but by the time of coming back to the nest, its spermatheca is filled with accumulated sperms, either fully or partially. After a mating flight, the queen bee turns onto fertilize the eggs by retrieving the sperms accumulated in the spermatheca.

### 3.3.2 Adoption of Biological Behavior of Honey-Bees for the Computational Purpose

For the execution of HBMO algorithm, the number of queens and drones need to be fixed in the beginning. The parameters associated with the queens that require pre-specification are maximum number of broods that can be produced after a mating flight, spermatheca size, initial speed, speed reduction factor and the maximum number of mating flights. Each honey bee in a colony represents a trial in the solution space. This point onwards, the words honey bee and solution are used synonymously. To begin the search, the ancestors of the honey-bee colony need to be initialized. The individuals are represented by the string of genes and the string length is equal to the number of decision variables. The genes of the individuals can be represented by binary coded or real coded values. To initialize the gene values, it can be assigned at random from the available discrete set of choices available in the list.

After initialization, the fitness of the individuals has to be evaluated. For minimization problems, the inverse of the objective or cost function value can be taken as the fitness of a honey bee. As analogous to the real honey-bee colony, the solution with more fitness value would act as queen.

The queen sets off the mating flight with an empty spermatheca and with a high initial speed. The mating flight can be treated as a set of transitions in a state space - the environment - where the queen moves between the different states in some speed and mates with the randomly selected drones one by one. The mating can take place only when the probabilistic rule of mating described in the equation (10) gets satisfied.

$$P = e^{-\frac{\Delta(f)}{S(t)}} \quad \ldots(10)$$

where,
$P$ = probability of mating;
$\Delta(f)$ = absolute difference between fitness values of the queen bee and the selected drone; and
$S(t)$ = speed of the queen.

After each transition in the space, speed of the queen gets reduced as per Equation (11).

$$S(t) = \alpha\, S(t-1) \quad \ldots(11)$$

where,
$\alpha$ = speed reduction factor.

To replicate the death of a drone after insemination, the solution which represents the drone that successfully mates with the queen bee has to be removed from the drone population. The mating process of a queen can hold on until the spermatheca gets filled with sperms or the speed falls below the threshold value. The mating flight is postulated for all the queens if there exists more than one queen.

After the mating flight, the new broods can be brought forth by coupling the queen's and drone's genes. In HBMO, the function of nurse bees is limited to brood care only as the new broods arise only from the queen's eggs. After bringing forth the specified numbers of broods, the objective function requires to be evaluated to subsequently compute the fitness values of the broods. The broods with higher fitness values than current queens' fitness, have to be replaced as queens of the next generation. The routine of mating flight and subroutine of brood generation need to be continued until the termination criteria is met. The termination criteria may be either arrival of maximum number of mating flights or no more improvement in solution over certain number of mating flights. The pseudo-code of this algorithm is presented in Figure 1.

### 3.3.3 Individual Representation
The number of genes in an individual is taken as the number of decision variables involved in the problem. That is, the number of pipes, whose diameters need determining, can be taken as the number of decision variables for a particular problem. For demonstration, the set of genes of a queen bee for an eight-variable problem can be represented as follows:

| Q1 | Q2 | Q3 | Q4 | Q5 | Q6 | Q7 | Q8 |

where Q1 to Q8 are the genes of the queen. The gene values of the individuals have to be initialized from the list of discrete variables at random.

The set of genes of a drone bee for an eight-variable problem can be represented as follows:

| D1 | D2 | D3 | D4 | D5 | D6 | D7 | D8 |

where, D1 to D8 are the genes of a drone.

```
1. initialization:
    define the number of queens and drones.
    define the number of genes of the individuals.
    define the maximum number of mating flights.
    define the spermatheca size of the queens.
    define the queen's initial speed and speed reduction factor.
    define the number of broods that can be produced by a queen after a mating flight.

2. generation of ancestors (initial colony)
    initialize the queens and drones gene values from the list of discrete variables at random.
    classification of honey-bees:
    evaluate the objective function and subsequently the fitness of the individuals.
    based on the fitness value, arrange the individuals in an ascending order.
    n=1;
    while n less than the number of queens
        individuals are queens.
        n=n+1.
    end while
    while n less than the number queens+drones
        individuals are drones.
        n=n+1.
    end while

3. generation of new colony:
    while maximum number of mating flight or no improvement in solution observed

            mating Flight:
            for each queen in the queen list
                while the queen's spermatheca has space or speed>min value
                    queen moves between states and randomly chooses drones
                    if a selected drone satisfies the probabilistic rule of mating,
                        add its sperm to the queen's spermatheca
                        remove the selected drone from the drone list
                    end If
                    update the queen's speed
                end While
            end for each queen

            fertilization:
            for each queen in the queen list
                while the queen's spermatheca empty
                    randomly retrieve the sperm from the spermatheca using uniform random
                    number generator
                        generate a brood by crossover the queen's and drone's genotypes
                        remove the used sperm from the spermatheca
                end while
            end for each queen

            update the queens:
            while the best brood is better than the worst queen
                replace the least–fittest queen with the best brood
                remove the best brood from the brood list
            end While
            update of drones list
    end while
```

**Figure 1 - Pseudo-Code of the HBMO Algorithm [9]**

The set of genes of a brood bee for an eight-variable problem can be represented as follows:

**B1  B2  B3  B4  B5  B6  B7  B8**

where B1 to B8 are the genes of a brood.

### 3.3.4 Crossover Operation - Single-Point Crossover

A single cutoff point has to be selected at random. For the brood's string, drone's genes have to be placed up to the cutoff point, after that the queen's genes have to be placed.

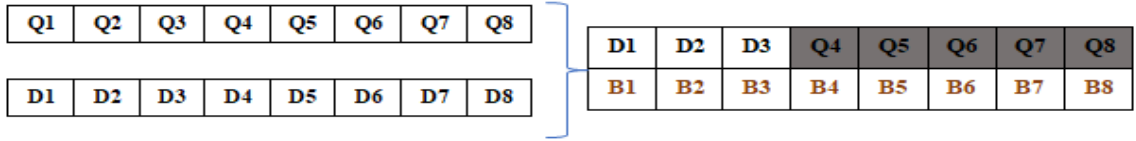

**Figure 2 - Random Cut-off Point Selection in the Queen's and the Drone's Genes [9]**

The random cut-off point selection is illustrated in Figure 2. As shown in Figure 2, B1, B2, B3, B4, B5, B6, B7 and B8 are the genes of a brood. The shaded boxes show the part of the string formed by the queen's genes.

### 3.4 Adapting the HBMO for the Designing of WDN

The procedure in which the HBMO was used for the optimal WDN design is presented in Figure 3. In order to begin the implementation of the algorithm, the WDN parameters, such as source head, nodal elevations and demands, minimum hydraulic-head value to be maintained at the demand nodes, maximum allowable friction loss gradient, pipe layout, length of the individual pipes, commercially available pipe sizes, and unit cost associated with the commercial pipes are fed to the model. After specifying these parameters, the initial colony of honey-bees, also known as set of ancestors, is generated. To evolve the ancestors, the gene values are assigned from the list of commercially available pipe sizes. After initialization, the queen is selected at random from the initial colony. A penalty is added with the objective function value of the solutions that do not satisfy the minimum hydraulic-head constraint. Similarly, another penalty is added with the objective function value of the solutions that do not satisfy the maximum friction loss gradient. The penalty values added for the constraint violation are modeled through equations:

$$(NP) = Max[(RH_j^{min} - RH_j), 0] * (NPF) \quad \ldots(12)$$
$$j = 1,2,3, \ldots, nd$$

$$(PP) = Max[0, (hf_i - hf_{max})] * (PPF) \quad \ldots(13)$$
$$i = 1,2,3, \ldots, np$$

where, NP, NPF, PP, and PPF denote nodal penalty, nodal penalty factor, pipe penalty, and pipe penalty factor, respectively. The solutions which do not satisfy the minimum hydraulic-head constraint and maximum friction loss constraint, are penalized in proportion to the hydraulic-head deficit $(RH_j^{min} - RH_j)$ and excess friction loss $(hf_i - hf_{max})$ values, as mentioned in Equation 12 and Equation 13, respectively.

The solutions which do not satisfy the minimum hydraulic-head constraint and maximum friction loss constraint, are penalized in proportion to the hydraulic-head deficit $(RH_j^{min} - RH_j)$ and excess friction loss $(hf_i - hf_{max})$ values, as mentioned in Equation 12 and Equation 13, respectively. The inverse of objective function value is taken as the fitness value of a solution. The individual with most fitness value is chosen as the queen and the honey-bees with lower fitness than the queen are then considered as drones.

Once the initial colony of honey bees are generated, the mating flight of the queen-bee begins. During the mating flight, the queen mates with the drones that are randomly selected from the drone population only when the probabilistic rule of mating gets satisfied.

Then, the queen starts to fertilize the eggs through crossover process. The fitness values of some of the newly born broods might be higher than that of the current queen. If a better brood appears, replacement of the queen by the brood with better fitness is carried out. The processes of mating flight, new brood generation, and replacement of the queen are continued until either the maximum number of mating flights or no improvement in solution over the pre-specified number of mating flights can be observed.

### 3.5 Gurudeniya Service Zone, Gurudeniya, Sri Lanka

The Honey-bee Mating Optimization (HBMO) algorithm was tested using the water distribution network of Gurudeniya Service Zone, Gurudeniya, Sri Lanka.

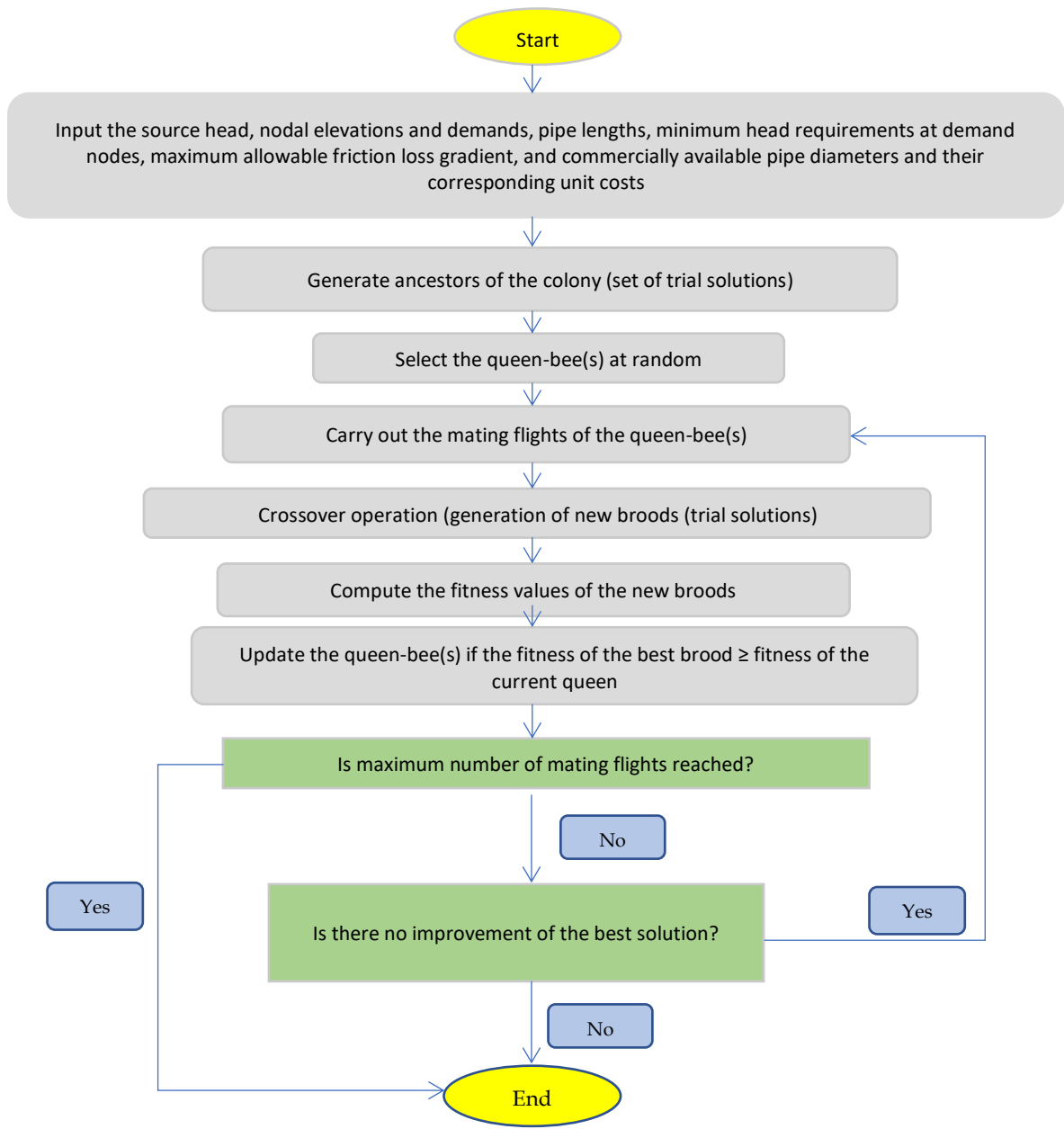

**Figure 3 - Flowchart for the Design of Water Distribution Network using Honey-Bee Mating Optimization (HBMO) Algorithm [9]**

The data set consisting of pipe layout, pipe lengths, nodal water demands, nodal elevations and reservoir (source) elevation data were obtained from National Water Supply and Drainage Board (NWS&DB), Sri Lanka. The schematic layout of the main pipe trunk, that has a serial pipe configuration, has been shown in Figure 4.

As shown in Figure 4, P1,P2,P3,…,P10 are the supply pipes. N1, N2, N3,…, N10 show consumer nodes and R1 shows the reservoir or source node. Lengths of pipes of the WDN are shown in Table 1. Nodal water demands and nodal elevations are shown in Table 2.

The list of commercially available pipe sizes and their unit costs (adapted from Cunha and Sousa [4]) are shown in Table 3. The elevation at the reservoir node has been taken to be 555m. Depending on the material of the pipes present in the WDN, Hazen William's coefficient was selected to be 130 for all the pipes.

**3.6   Assignment of Values for the Parameters in the Algorithm**

The initial step when implementing HBMO was to evolve the ancestors. For this, the genes of individuals were assigned from the randomly selected pipe sizes.

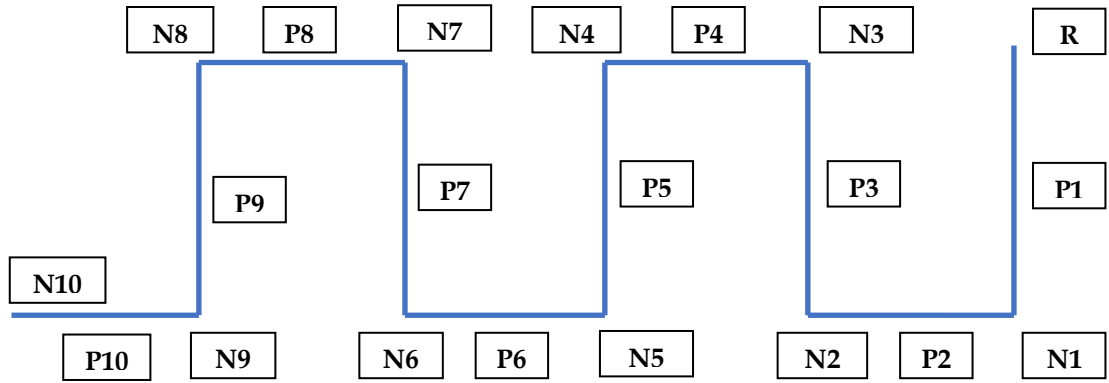

**Figure 4 - Schematic Layout of Main Pipe Trunk – Gurudeniya Service Zone, Gurudeniya, Sri Lanka**

**Table 1 - Lengths of Pipes of the Gurudeniya WDN**

| Pipe ID | Length / (m) |
|---|---|
| P1 | 690 |
| P2 | 1120 |
| P3 | 120 |
| P4 | 270 |
| P5 | 630 |
| P6 | 280 |
| P7 | 420 |
| P8 | 230 |
| P9 | 290 |
| P10 | 980 |

**Table 2 - Nodal Water Demands and Nodal Elevations of Gurudeniya WDN**

| Node ID | Demand / (m³/day) | Elevation / (m) |
|---|---|---|
| N1 | 796.52 | 452 |
| N2 | 127.50 | 517 |
| N3 | 112.50 | 519 |
| N4 | 165.00 | 535 |
| N5 | 258.76 | 490 |
| N6 | 131.25 | 481 |
| N7 | 168.76 | 476 |
| N8 | 228.76 | 486 |
| N9 | 333.76 | 462 |
| N10 | 37.50 | 480 |

**Table 3 - Commercially Available Pipe Sizes and their Unit-Costs [4]**

| Diameter/ (mm) | Unit-Cost /(Units) |
|---|---|
| 25.40 | 2 |
| 50.80 | 5 |
| 76.20 | 8 |
| 101.60 | 11 |
| 152.40 | 16 |
| 203.20 | 23 |
| 254.00 | 32 |

The number of genes of an individual was made equal to the number of pipes in the system. Therefore, for Gurudeniya WDN, each honey-bee was provided with ten number of genes for each. The analysis was carried out for a combination of one queen-bee and 499 drones. The number of mating flights was selected to be 1000. The spermatheca size of the queen-bee was selected as 100. The initial speed and the speed reduction factor were taken as 2 and 0.95, respectively. The threshold probability of the queen-bee to mate with a drone was selected as 0.01. In this state of the study, the Pipe Penalty Factor was kept constant, and only the Nodal Penalty Factor was varied to check its effect on the solution. At each of the nodes, the minimum residual head ($H_R^{min}$) required was set to 10 m, while the maximum friction loss ($h_{f\ max}$) allowed was set to 0.005 m/m, as per the design standard adapted at NWS&DB, Sri Lanka.

## 4. Results and Discussion

Solutions that were obtained for six trials are shown in Table 4. It is clear that, all six trials show closeness to the solution implemented by NWS&DB, Sri Lanka. It can be intuitively stated that, sudden elevation rises in the nodes N2, N4, N8, and N10 may have influenced diameters of pipes P2, P4, P8, and P10 to get noticeable fluctuating values from Trial 1 to Trial 6.

Table 5 shows the total penalty values and total cost values for different values of Nodal Penalty Factors assigned in the algorithm. Here, the Total Penalty is the summation of Nodal Penalty and the Pipe Penalty described in Equations 12 and 13, respectively. It could be noticed that, for small Nodal Penalty Factor values in Trial 1 and Trial 2, the solutions given by the algorithm are violating the hydraulic constraints of the problem by indicating large Total Penalties. Upon increasing the Nodal Penalty Factor value, zero Total Penalties have been indicated in Trials 3, 4, 5, and 6. With the extra increase of the Nodal Penalty Factor, it can be observed from the Trial 6 that, Total Cost has risen to a relatively high value. It can be undeniably observed that, when the Nodal Penalty Factor is in close proximity to 5000, the cost of solution by the algorithm has gone below the cost of the solution NWS&DB has used in project implementation. Therefore, it can be promising that, with the use of HBMO algorithm, more solutions could be obtained whose costs are below the cost value of NWS&DB implemented solution.

## 5. Conclusions

In this study, the use of Honey-Bee Mating Optimization (HBMO) for the design of WDN in the Gurudeniya Service Zone, Gurudeniya, Sri Lanka, has been examined. From the implementation of the HBMO algorithm, it could be observed that the HBMO algorithm is successful in addressing the WDN design problem with given mathematical constraints. In-built hydraulic constraints have enabled the algorithm to deliver results without needing a separate hydraulic simulation software in the design. The comparison of results given by HBMO algorithm to the solution implemented by NWS&DB shows that the HBMO algorithm not only gives feasible solutions, but also gives solutions that have lower costs involved when compared with NWS&DB implemented solution. It should be noted that, although this study has used the unit costs involved with the pipes, the true costs could be significantly high, depending on the value of a monetary unit concerned.

Table 4 - Solutions Obtained through Six Numbers of Trials and the Solution that was Implemented by NWS&DB (All Values are in mm)

|     | Trial 1 | Trial 2 | Trial 3 | Trial 4 | Trial 5 | Trial 6 | NWS&DB Implemented Solution |
|-----|---------|---------|---------|---------|---------|---------|------------------------------|
| D1  | 203.2   | 203.2   | 254.0   | 203.2   | 254.0   | 203.2   | 254.0                        |
| D2  | 203.2   | 203.2   | 203.2   | 203.2   | 203.2   | 254.0   | 203.2                        |
| D3  | 101.6   | 101.6   | 152.4   | 254.0   | 203.2   | 152.4   | 203.2                        |
| D4  | 254.4   | 152.4   | 152.4   | 203.2   | 203.2   | 203.2   | 152.4                        |
| D5  | 101.6   | 152.4   | 101.6   | 152.4   | 152.4   | 203.2   | 203.2                        |
| D6  | 203.2   | 152.4   | 254.0   | 152.4   | 101.6   | 101.6   | 101.6                        |
| D7  | 76.2    | 101.6   | 152.4   | 203.2   | 101.6   | 203.2   | 101.6                        |
| D8  | 203.2   | 152.4   | 101.6   | 152.4   | 76.2    | 254.0   | 76.2                         |
| D9  | 152.4   | 203.2   | 101.6   | 203.2   | 101.6   | 152.4   | 76.2                         |
| D10 | 76.21   | 76.21   | 152.4   | 25.4    | 50.8    | 76.2    | 76.2                         |

Table 5 - Total Penalty and Total Cost against Nodal Penalty Factor for the Six Trials

| Trial Number    | Nodal Penalty Factor | Total Penalty | Total Cost / (Units) |
|-----------------|----------------------|---------------|----------------------|
| 1               | 2000                 | 11846         | 86090                |
| 2               | 3000                 | 10586         | 84640                |
| 3               | 3900                 | 0             | 98090                |
| 4               | 4000                 | 0             | 88210                |
| 5               | 5000                 | 0             | 84520                |
| 6               | 6000                 | 0             | 95830                |
| NWS&DB solution |                      |               | 89110                |

Moreover, having proven the applicability of the HBMO algorithm for serial pipe configuration in Gurudeniya Service Zone, authors recommend further research on checking the applicability of the same algorithm on looped pipe configuration.

## Acknowledgements

Authors would like to thank the National Water Supply and Drainage Board for their data being used in this study.